\documentclass[preprint]{elsarticle}

\begin{document}

\title{Averaging Transformations of Synaptic Potentials on Networks}

\author[rvt]{H. R. ~Noori\corref{cor1}}
\ead{hamid.reza.noori@iwr.uni-heidelberg.de, hnoori@princeton.edu}
\cortext[cor1]{Corresponding author}
\address[rvt]{Interdisciplinary Center for Scientific Computing, University of Heidelberg, Im Neuenheimer Feld 294, 69120 Heidelberg, Germany}

\begin{abstract}
The problem of the transformation of microscopic information to the macroscopic level is an intriguing challenge in computational neuroscience, but also of general mathematical importance. Here, a phenomenological mathematical model is introduced that simulates the internal information processing of brain compartments. Synaptic potentials are integrated over small number of realistically coupled neurons to obtain macroscopic quantities. The striatal complex, an important part of the basal ganglia circuit in the brain for regulating motor activity, has been investigated as an example for the validation of the model.
\end{abstract}

\maketitle

\section{Introduction}

The brain nuclei, as parts of complex brain networks, are comprised by different types of inter- and projection neurons within subnetworks of highly complex structure. 
Beside the anatomical formation, functional properties in integration and processing of neural information become adapted during the development of these regions until adolescence. In this study, adult brain nuclei are considered after their final development, as part of information processing pathways. 

Experiments on the topology of brain regions suggest that these structures can be categorized as heterogeneous media. There are various mathematical methods describing dynamical processes in heterogeneous media, such as asymptotic analysis, homogenization, and integral equations (Pavliotis and Stuart \cite{Pa}, Haken \cite{Ha}, E and Engquist \cite{E}, Levin and Chao\cite{Le}). There are traditional network approaches on the effects of properties such as the local connectivity of neurons on the striatal function (Wickens et al. \cite{Wi}). These studies concentrate merely on the network structure of a particular brain region e.g. striatum and even then do not contain some important features such as the role of large spiny neurons in the striatal function. Large spiny neurons in the striatum are cholinergic and are of great importance for processes involving addiction and food intake (Rada et al. \cite{Rad}, \cite{Rad1}, Avena et al. \cite{Ave}).

In the present study, simplified but essential physiological processes have been chosen as the biological foundation for the mathematical model. The reader is advised to consider the appendix for a brief introduction into the biological terminology. 

The course of information processing in brain regions depends on the propagation of synaptic potentials - which are averaged postsynaptic potentials - along the subnetworks of these regions. Thus, understanding the physiology of this integration process requires the comprehension of the neuronal architecture of the brain region of interest. The regional subnetworks consist of different types of neurons and characteristic connections among them that can be obtained by experiments. In this study the intraneuronal connections are subdivided into the global, and ultrastructural morphology. Golgi and immunohistochemical methods provide optimal frameworks for obtaining the values of the two morphological dimensions experimentally. The immunohistochemical studies reveal the appearance of certain neurotransmitter systems in a structure; Golgi studies provide information about the form and type of neurons expressing the obtained neurotransmitter systems and on the morphological dimensions of the brain region in terms of the characteristic connections between different types of neurons in a brain region. In general, the structural topology of the brain region is then obtained. The knowledge about the appearance of different types of neurons and the intraneuronal synaptic interactions inside a nucleus allow us to abstract its structure through networks. 

One way of representing such networks is through their realization in continuous spaces. The attendance of a high amount of neurons supports the idea of the appropriateness of such an approach for this purpose. 

Because of the importance of the global and ultrastructural morphology for the information integration process in a nucleus and the dependency of their values on the spatial distribution of neurons, they shall be modelled as spatial variables parametrized by the neuron (neurotransmitter) type. 

To provide a general framework for the regional neuronal activity, we introduce an integral operator averaging the information on two refinement scales. The idea is to cover the space including the nucleus with copies of discrete fundamental domains consisting of neuronal assemblies which we call n-cells. Roughly spoken, n-cell patches are finite networks of neurons. The ultrastructural morphology of a nucleus in terms of the synaptic connectivities between different neurons provides information on the local connectivities in the n-cells, which are embedded as edges of neuronal network inside the n-cells. Within the information on the adherence functions of n-cells (global morphology variables) and the distribution of neurons, the averaging of synaptic potentials across the network of nuclei is completed. The main advantage of this method is in its appropriateness for applications on various brain regions.  
     
\begin{figure}
\begin{center}
\includegraphics[width=80mm]{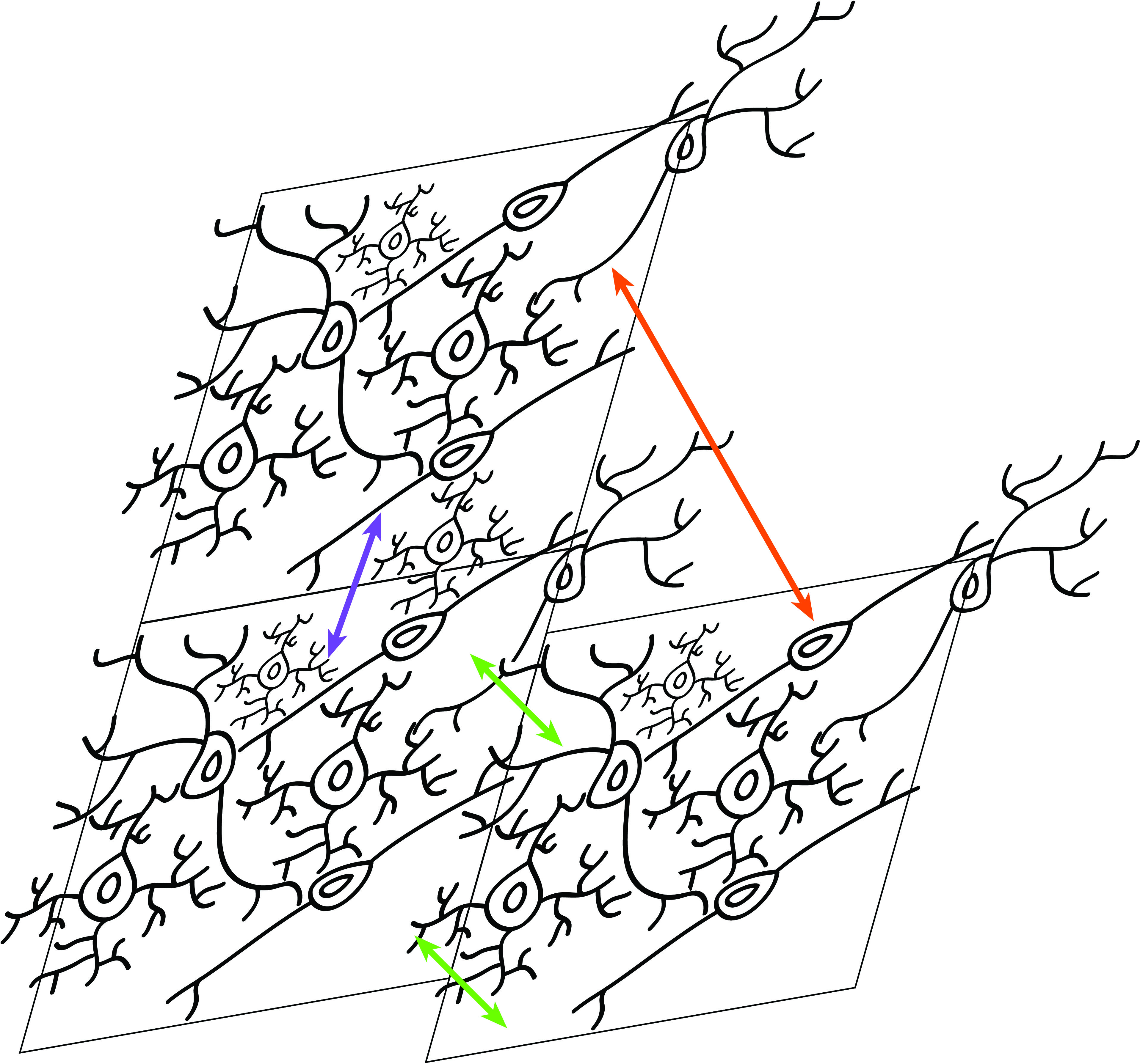}
\vspace{2.5cm}
\caption{Composition of n-cells}
\label{single}
\end{center}
\end{figure}

Assembling a suitable number of n-cells by considering the large-scale interaction between the n-cells (adherence of n-cells), we obtain information on the combined influence of several synaptic inputs on macroscopic neuronal activity (Fig. \ref{single}). One trivial example is the subthalamic nucleus which consists of only one class of interneurons. In this case the number of n-cells under consideration is one (Parent and Parent \cite{Par}).    

The aim of this study is to present a mathematical technique for averaging microscopic information along complex network structures such as the networks of single brain regions. The technique is applied on the corpus striatum - a part of the basal ganglia circuit which plays an important role in the regulation of motor behaviour- and shown to mimic faithfully the qualitative oscillatory behaviour of this brain region. The general formulation of the averaging operator suggest its suitability for other discrete multiscale problems, especially in the research area of material sciences.

\section{Model Statement} 

\paragraph{Statement:} The integration of the electric activity in a brain nucleus (information processing) depends on the distribution of its comprising neurons, the ultrastructural morphology as a local variable quantity and the global quantity in terms of the composition of neuronal fundamental domains (n-cells), that describes the topology of the brain region.

The concept of the n-cells is very essential for the understanding of the averaging transformation. The idea is to path up a region with multiple copies of discrete fundamental domains called n-cells. The n-cells compose the whole brain region of interest by applying the global morphological properties and contain the local morphological properties as finite graphs. Hence, they unify the morphological dimensions into one notion. 

\bigskip\noindent

\textbf{Definition 1:}
A n-cell is a neuronal assembly represented by a finite graph characterized by the local ultrastructural morphology. Within proper composition of the n-cells, as discrete fundamental domains, they patch up the region of interest. 

\bigskip\noindent

Let $\Omega \subset \mathbf{R}^3$ a compact, path-connected set, be the brain region of interest. Let $\Omega_i \subset \Omega$, $i \in I$ be a skeleton representation of the finite graph of a n-cell, and $I$ be the finite index set with $\#I = \# \{\mbox{n-cells} \}$. 
Furthermore, let $x \in \Omega_i$ be the variable describing the network position of a neuron relative to the n-cell $\Omega_i$; and $y$ the continuity variable describing the spatial position of a neuron embedded as points in $\Omega$ as a subset of the complete space $\mathbf{R}^3$. $\{z\}$ denotes the finite set of neurotransmitter families that appear in the brain region of interest. By allocating any neuron to a neurotransmitter, $z$ parametrizes the proposed averaging. 

\bigskip\noindent

\textbf{Definition 2:}
A synaptic potential $u_x^z(t)$ is the sum of all Hodgkin-Huxley postsynaptic potentials of a neuron of type $z$ at position $x$ inside a n-cell. 
$$\forall i \in I \ , \ \forall z \ , \ \forall x \in \Omega_i \ : \ u_x^z: [0,T] \rightarrow \mathbf{R} \ .$$

\bigskip\noindent

To average the synaptic potentials $u_x^z(t)$ in n-cells, the local connectivity function $\psi^z : \Omega_i \rightarrow \mathbf{R}^+$ has to be characterized. By introduction of the global connectivity function $\chi_i: \Omega \rightarrow \mathbf{R}^+$, describing the composition of n-cells, we will cover $\Omega$ will n-cell patches. We write 
$\chi:=(\chi_i)$ for the $|I|$-vector of composition mappings.
The normalized distribution function $\rho^z: \Omega \rightarrow [0,1]$ as the kernel of the averaging operator will complete the information required for the integration of synaptic potentials across the network of brain regions. 
These functions $(\psi^z, \chi_i, \rho^z)$ are obtained from Golgi-  and immunohistochemical studies.

The continuous realization of the network structure of the brain nuclei and the idea of transforming the microscopic, synaptic information to the macroscopic level compatible with the internal structure of a brain region suggest the construction of an averaging integral operator along networks. 

\bigskip\noindent

\textbf{Definition 3:}
A compartment is a multiple $C:=(\Omega, \rho^z, \chi, \psi^z, I)$ which characterizes a brain nucleus by its ultrastructural and global morphology. 

\bigskip\noindent

\textbf{Definition 4:}
Let $C$ be a compartment, then the averaged potential $v(t)$ of synaptic potentials $u_x^z(t)$ across the network of a single brain nucleus is represented by:

\begin{equation}
v(t) \ = \ \sum_{z} \int_{\Omega}{\rho^z(y) \cdot \frac{\sum_{i \in I}{\chi_i(y) (\sum_{x \in \Omega_i}{\psi^z(x) u^z_x(t)})}}{\left| I \right| g(\chi(y))} dy} \ , \ \bigcup_i \Omega_i \ \subseteq \Omega \ ,
\end{equation}
$g: \mathbf{R}^{|I|} \rightarrow \mathbf{R}^+$ is the proper averaging parameter estimated by immunohistochemical experiments. It provides a combinatorial dimension for the composition of the n-cells in a compartment.

\bigskip\noindent   

The universality of the concept of a compartment provides the possibility of applying the averaging transformation of this type for a wide category of network multiscale problems.

\section{Numerical Simulations}

The aim of this section is the numerical investigation of the dynamical behaviour of a compartment by known 
synaptic dynamics (in terms of electric activity $u_x^z(t)$) of single neurons. The corpus striatum is used as an non-trivial example to illustrate the efficiency of the averaging operator method. The complexity of the striatal structure on the one hand, and the importance of this brain region as a substrate of the basal ganglia in the regulation of motor activity and neurological diseases on the other hand, make it a proper representative for such investigations. First,  the main morphological and ultrastructural properties of this compartment are summarized which are required to comprise the integral function. These properties include the spatial distribution of the different neural populations ($\rho^z$), the neurochemical classification ($z$), and the intraneuronal connections ($(\psi^z, \chi_i, g)$). Then, the simulation results are discussed and compared with the experiments. 

A category of proper experiments for the validation of the numerical results is represented by the local field potential (LFP) studies. Thereby, a signal is recorded using a low impedance extracellular microelectrode, placed sufficiently far from individual local neurons to prevent any particular cell from dominating the electrophysiological signal. This signal is then low-pass filtered, cut off at $\approx 300 Hz$, to obtain the local field potential (LFP). The low impedance and positioning of the electrode allows the activity of a large number of neurons to contribute to the signal. The unfiltered signal reflects the sum of action potentials from cells within approximately $50-350 \mu m$ from the tip of the electrode (Legatt et al. \cite{Legat}). 

The frequency of the LFP signals in corpus striatum will then be compared with those appeared as the results of the numerical simulation of the averaging transformation by given morphological properties.
     
\subsection{The Corpus Striatum}

\begin{figure}
\begin{center}
\includegraphics[width=80mm]{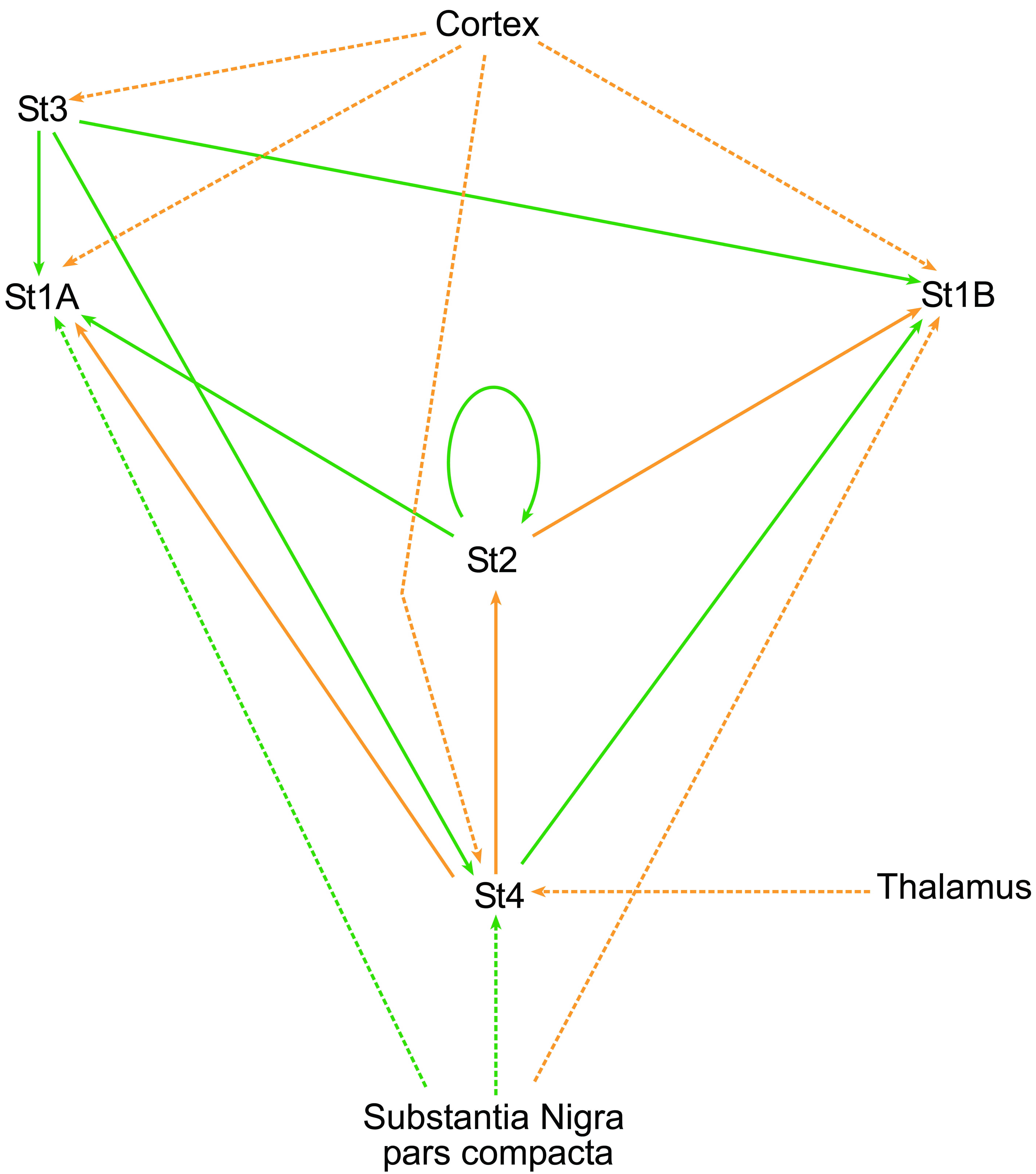}
\vspace{2.5cm}
\caption{Schematic ultrastructural morphology of the primate's striatum, including the neurochemical synaptology. The green/orange arrows denote inhibitory/excitatory afferents. The dashed arrows are external efferents to the striatal neurons. Hereby, $St_{1A}$ and $St_{1B}$ denote the spiny projection neurons which project to Globus Pallidus externa and Globus Pallidus interna. $St_2$ denotes the dopaminergic interneurons, $St_3$ the GABAergic interneurons, and $St_4$ the acetylcholinergic neurons.}
\label{striatal}
\end{center}
\end{figure}

The neuronal populations of the striatum could be divided into four classes: Spiny projection neurons (about $96 \%$ of the whole neural population) are GABAergic neurons, which get external inputs from cortical areas, thalamus and substantia nigra pars compacta (Smith et al. \cite{Smi}). These neurons also get inputs from the dopaminergic, GABAergic and acetylcholinergic interneurons. Large aspiny neurons (cholinergic), GABAergic interneurons and dopaminergic neurons comprise the rest population of striatal neurons. The cholinergic interneurons get external inputs from thalamus and substantia nigra pars compacta; and internal inputs from the GABAergic interneurons. They project to the dopaminergic and projection neurons. The action of the cholinergic interneurons is assumed to have opposite effects on the action of the dopaminergic neurons on the whole network. The cholinergic neurons inhibit the activity of the projection neurons of the direct pathway (mostly $D_1$ system) and disinhibit the activity of the indirect pathway neurons ($D_2$ system) (DiFiglia \cite{Di2}, Graybiel et al. \cite{Gra}, DiFiglia et al. \cite{Di0}, DiFiglia and Carey \cite{Di1}, Gerfen \cite{Ge}, Betarbet et al. \cite{Bet}, Flores-Hernandez et al. \cite{Fl}, Graveland et al. \cite{Gr}, Bennett and Wilson \cite{Be}). A schematic representation of the microcircuitry of corpus striatum is shown in Fig. \ref{striatal}. 

MATLAB has been used to simulate the averaging across the network of striatum. An overall number of $6400$ neurons has been used for the simulation of averaged potentials in striatum. The density distribution function $\rho$ is represented in Figure 3. The local and global connectivity function $\psi$ and $\chi$ defined the necessary computational parameters based on Figure 2. We have induced a single excitatory signal from cortical projection neurons to the cholinergic neurons in striatum to simulate the simplest integration process. 

\begin{figure}
\begin{center}
\includegraphics[width=\textwidth]{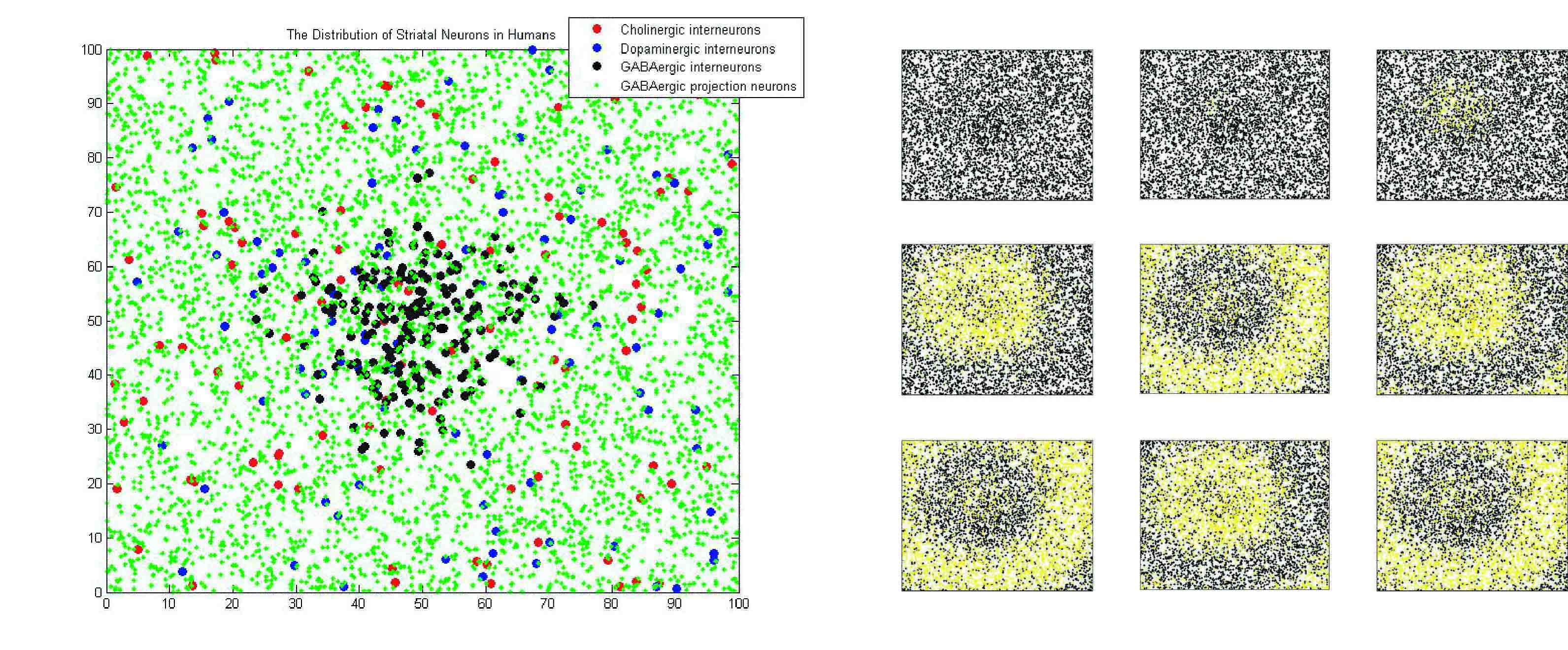}
\vspace{2.5cm}
\caption{Left: The distribution of striatal interneurons based on experimental data. Right: The oscillatory dynamics of the activity of striatal neurons by a single constant activation with a realistic spatial distribution of the neurons (tonically active cholinergic neurons). The light nodes denote active neurons.}
\label{striatal2}
\end{center}
\end{figure}

It reveals that the activity transmission induced by a constantly activated neurons is approximately radial and oscillatory. The activity of a single neuron, forced by a constant input, produces activation-waves on the striatal populations (Fig. \ref{striatal2}). 

Such neurons are simplified versions of the tonically active cholinergic interneurons. Activations of neurons of other classes suggest similar activation waves across the striatal populations. We observe that the simulated oscillations of $v(t)$ ($\approx 50 Hz$) are qualitatively correlated with the LFP-studies of basal ganglia which also suggest oscillations in control patients (Boraud et al. \cite{Bo}). The efficiency of the averaging transformation has been also investigated for other brain compartments. In general, the averaged potentials of the compartments reveal oscillatory dynamical behaviour. For example, synchronized gamma oscillations were observed in the globus pallidus. The simulation results for the globus pallidus could be well reproduced by the reader using the averaging transformation. Therefore, the representation of these results is omitted in the present study.

\section{Discussion}

In conformance with LFP experiments, the multiscale averaging operator (1) transforms the synaptic potentials along the network structure of brain region such as corpus striatum to local field potentials of the same frequency range. This represents a first step from the neuroscience at Hodgkin-Huxely dimension to the system biological level of consideration. 
Although this work was initially inspired by one of the most intriguing problems of computational neuroscience, it appears to be appropriate for several other intriguing applications. The characterization of the discrete fundamental domains in a complex network in terms of n-cells which is one of the dominating concepts of this mathematical approach, provides the possibility of the multiscale transformations of information along discrete structure that are continuously realized. 
Beside the advantages of this model such as its agreement with electrophysiological experiments, further computational investigations are needed to verify its validity for more complex systems. 

\bigskip\noindent

\paragraph{Acknowledgment.}
The financial support by the Interdisciplinary Center for Scientific Computing and the International Graduiertenkolleg 710 (DFG) and the funding by National Genome Research Network are acknowledged.  

\bigskip\noindent

\paragraph{Appendix.}

\bigskip\noindent
\begin{itemize}
\item \textit{Neurotransmitter}-Neurotransmitters are the most common class of chemical messengers in the nervous system;
\item \textit{Synapse}-Chemical synapses are specialized junctions through which neurons signal to each other and to non-neuronal cells such as those in muscles. Chemical synapses allow neurons to form circuits within the central nervous system. They are crucial to the biological computations that underlie perception and thought;
\item \textit{Morphology}-The term morphology in biology refers to form, structure and configuration of an organism. This includes aspects of the outward appearance (shape, structure, colour, pattern) as well as the form and structure of the internal parts like bones and organs;
\item \textit{Ultrastructural morphology}-Cell and tissue morphology in the electron-microscope level;
\item \textit{Interneuron}-An interneuron (also called relay neuron, association neuron or local circuit neuron) is a multipolar neuron which connects afferent neurons and efferent neurons in neural pathways. Like motor neurons, interneuron cell bodies are always located in the central nervous system;
\item \textit{Golgi method}-Golgi's method is a nervous tissue staining technique;
\item \textit{Immunohistochemical studies}-Immunohistochemistry or IHC refers to the process of localizing proteins in cells of a tissue section exploiting the principle of antibodies binding specifically to antigens in biological tissues;
\end{itemize}

\end{document}